\documentclass[10pt,a4paper]{amsart} 
\usepackage{graphicx,latexsym,amssymb,color}
\theoremstyle{plain} 
\newtheorem{thm}{Theorem}[section]

\newtheorem{prop}[thm]{Proposition}

\newtheorem{Problem}[thm]{Problem}
\theoremstyle{definition}
\newtheorem{defi}{Definition} 
\theoremstyle{remark} 
\newtheorem{rem}[thm]{Remark}

\numberwithin{equation}{section}
\numberwithin{figure}{section}
\newcommand{\bd}{\begin{description}}   
\newcommand{\ed}{\end{description}} 
\newcommand{\bF}{\begin{figure}[!h]}      \newcommand{\eF}{\end{figure}} 
\newcommand{\ba}{\begin{array}}      \newcommand{\ea}{\end{array}} 
\newcommand{\bc}{\begin{center}}     \newcommand{\ec}{\end{center}} 
\newcommand{\be}{\begin{enumerate}}  \newcommand{\ee}{\end{enumerate}} 
\newcommand{\beq}{\begin{eqnarray}}  \newcommand{\eeq}{\end{eqnarray}} 
\newcommand{\beQ}{\begin{eqnarray*}} \newcommand{\eeQ}{\end{eqnarray*}} 
\newcommand{\bi}{\begin{itemize}}    \newcommand{\ei}{\end{itemize}}

\newcommand{\Sm}{\Sigma} 
\newcommand{\s}{\sigma}

\begin{document} 
\title[On the Classification of Links up to Finite Type]{On the classification of links up to finite type} 

\author[N. Habegger]{Nathan Habegger (Ju Li Wu Xiao-Song)}
\address{UMR 6629 du CNRS, Universit\'e de Nantes \\ D\'epartement de Math\'e\-matiques \\
         2 rue de la Houssini\`ere \\ 44072 NANTES Cedex 03, France}
         \email{habegger@math.univ-nantes.fr}

\author[J.B. Meilhan]{Jean-Baptiste Meilhan} 
\address{CTQM -- Department of Mathematical Sciences\\
         University of Aarhus \\
         Ny Munkegade, bldg. 1530\\
         8000 Aarhus C, Denmark}
	 \email{meilhan@imf.au.dk}
%
\begin{abstract} 
We use an action, of $2l$-component string links on $l$-component string  links, defined by Habegger and Lin, to lift the indeterminacy of finite type link invariants.  The set of links up to this new indeterminacy is in bijection with the orbit space of the restriction of this action to the stabilizer of the identity. \\
Structure theorems for the sets of links up to $C_n$-equivalence and Self-$C_n$-equivalence are also given.
\end{abstract}

\maketitle

\centerline{(In fond rememberance of Xiao-Song Lin, 1957-2007) }

\section{Introduction} \label{intro}
In \cite{M1,M2}, Milnor defined invariants of links, known as the Milnor $\overline \mu$ invariants.  In fact, these invariants are not universally defined, \emph{i.e.,} if the lower order invariants do not vanish, they are either not defined, or at best, they have indeterminacies.  

In \cite{HL1}, the notion of string link was introduced, together with the philosophy that Milnor's invariants are actually invariants of string links.  Indeterminacies are determined precisely by the indeterminacy of representing a link as the closure of a string link.  This philosophy led to
the classification of links up to homotopy, and to an algorithm constructed by Xiao-Song.  (Here and throughout, we will often refer to Xiao-Song Lin by his first name.)
More precisely, Xiao-Song and the first author constructed an orbit space structure for the set of links
up to homotopy.  The group action was `unipotent', meaning it acted trivially 
on the successive layers of the nilpotent homotopy string link group.  
This was the determining structural feature which
underlay the successful construction of Xiao-Song's algorithm.

In \cite{HL2}, an analogous orbit space structure for link
concordance was obtained and a study of the algebraic part of link
concordance, corresponding to the Milnor concordance invariants, was made. The theory also applies 
to more general `concordance-type' equivalence relations, in particular to those
studied by Kent Orr \cite{O} and developed in Xiao-Song's thesis.

With the advent of the physical interpretation of the Jones Polynomial \cite{J}, predicted by
Atiyah \cite{A} and established by Witten \cite{W}, 
a whole new area, known as Quantum Topology, emerged.  
Its perturbative aspects 
are succinctly summarized in the Universal Finite Type Invariant known as the 
Kontsevich Integral \cite{K}. 

Recall that, in the seminal paper \cite{L}, Xiao-Song had shown that Milnor Invariants are finite type invariants of string links.
We refer the reader to the paper of the first author and Gregor Masbaum 
\cite{HM}, where a formula is given which computes the Milnor Invariants directly from the Kontsevich Integral.

No successful attempt has been made at applying the methods of \cite{HL1,HL2} to the
Vassiliev invariants \cite{V}.  (The Vassiliev invariants were shown by Xiao-Song and Joan Birman \cite{BL} 
to be those invariants which satisfy the properties of
finite type invariants.  Subsequently, Bar-Natan \cite{B} adopted those properties as axioms for
finite type invariants.)
This is so, because, as we show, the classification scheme does \emph{not} hold.
Thus, in the philosophy of \cite{HL1,HL2}, the finite type invariants of links ought
to be refined.

In this paper, we make such a refinement and show that after refinement, 
the classification scheme applies.  We also show it applies to  $C_n$-equivalence
and to Self-$C_n$-equivalence.
 
\subsection*{Acknowledgments}
The authors extend thanks to J\o rgen Andersen and Bob Penner 
for organizing the stimulating conference \emph{Finite Type Invariants, Fat Graphs and Torelli-Johnson-Morita Theory} 
at the CTQM in Aarhus, March 2008, during which this work was done.
They also thank J\o rgen Andersen and Gw\'ena\"el Massuyeau for useful conversations.

The first author thanks the University of Nantes for releasing his time
for research. He also wishes to thank Zhenghan Wang for the invitation to speak 
at the conference in honor of Xiao-Song
at the Chern Institute in Tianjin, July 2007.  It was there that he was baptized
by Xiao-Song's wife, Jean, during a short ceremony in which he received his 
Chinese name, Ju Li Wu Xiao-Song.

\section{Preliminaries}

Let $D^2$ be the standard two-dimensional disk, and let $I$ denote the unit interval.  
Recall from \cite{HL1} the notion of string link.  
\begin{defi}
Let $l\ge 1$.  An $l$-component string link is a proper embedding,  
\[ \sigma : \bigsqcup_{i=1}^l I_i \rightarrow D^2\times I, \]
of the disjoint union $\coprod_{i=1}^{l} I_i$ of $l$ copies of $I$ in $D^2\times I$, such that the $j=0,1$ levels are preserved and
$\partial_j \s \subset D^2\times \{  j\}$ is the  standard inclusion of $l$ points in $D^2$.  
By an abuse of notation, we will also denote by  $\s \subset D^2\times I$,  the image of the map $\s$.  
\end{defi}
Note that we do not require that the $t$ levels, for $t\in I$, be preserved.  A string link is a pure braid precisely when
it preserves the $t$ levels for all $t\in I$.  
Note also that each string of an $l$-component string link is equipped with an (upward) orientation induced by the natural orientation of $I$.

\begin{figure}[!h]
\includegraphics{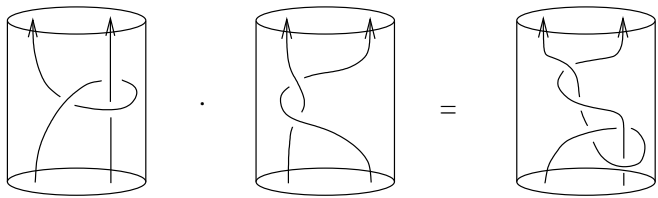}
\caption{Multiplying two $2$-component string links.}\label{SL}
\end{figure}

The set $SL(l)$ of isotopy classes of $l$-component string links (fixing the boundary) has a monoidal structure, with composition 
given by the \emph{stacking product} and with the trivial $l$-component string link $1_l$ as unit element.   See Figure \ref{SL}.  

\begin{rem}
In the above, one may replace the disk $D^2$ with any surface $S$ to get the notion of a string link in $S\times I$.  The
$l$-component string links in $S$, up to isotopy, again has a monoidal structure.

\end{rem}

We denote by  $L(l)$ the set of isotopy classes of $l$-component links.  
By a link, we mean an embedding $\coprod_{i=1}^{l} {\bf S}^1_i \rightarrow {\bf R}^3$.  Thus the
components are ordered and oriented.  There is an obvious surjective \emph{closure} map 
\[ \hat\ :  SL(l)\longrightarrow L(l) \]
which closes an $l$-component  string link $\s$ into an $l$-component link $\hat\s$.  

\begin{figure}[!h]
\includegraphics{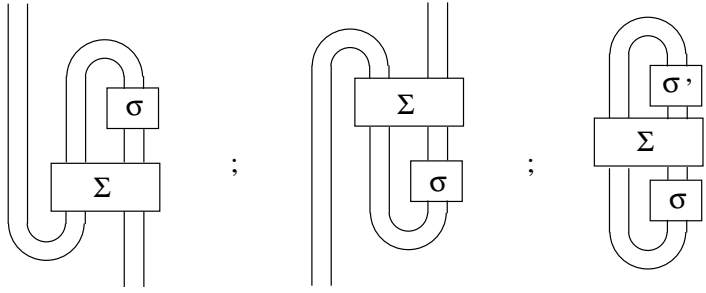}
\caption{Schematical representations of the left and right actions of $\Sigma$ on $\s$, 
$\Sigma \s$ and $\s \Sigma$,
and of the closed link $\s \Sigma \s'$.}\label{leftright}
\end{figure}

In \cite{HL1}, Xiao-Song and the first author introduced a certain left (resp. right) action of the monoid of isotopy classes of $2l$-component
string links on $l$-component string links.  
See Figure \ref{leftright} for an illustration of these actions.  
Thus given  two $l$-component string links $\s$, $\s'$, and a  $2l$-component string link $\Sm$, one has $l$-component string
links  $\Sm \s$, 
$\s\Sm $,
and a closed link  $\s\Sm \s'$.

One may represent the closure $\hat\s$ of a
string link $\s$ as $1_l 1_{2l}\s$, as well as $\s 1_{2l}1_l$, and also as $1_l (1_l\otimes\s) 1_l$, where $\s_1\otimes\s_2$ denotes the $2l$-component string link
obtained by horizontal juxtaposition.  (One orients all strands appropriately in the above, e.g., in $\Sigma$, one must reverse the parametrization of the first $l$ strands.) 

The following result on basing links was proven in \cite{HL2}.

\begin{prop}\label{basing} 
Let $\s_1$, $\s_2$ be two $l$-component string links whose closures are isotopic. Then there is a
$2l$-component string link $\Sm$, with  $1_l\Sm $ isotopic to $\s_1$, and $\Sm 1_l$  isotopic to $\s_2$. 
\end{prop} 

\section{The Habegger-Lin Classification Scheme}  
In \cite{HL2} a structure theorem was proven for certain `concordance-type' equivalence relations on the set
of links.   Given here for the convenience of the reader, 
though stated slightly differently, the result is in fact implicit in the proof in \cite{HL2}.

Consider an equivalence relation $E$ on string links and on
links (for all $l$), which is implied by isotopy.   We will denote by $E(x)$, the $E$ equivalence class of
$x$.  We denote by $ESL(l)$, resp. $EL(l)$,
the set of $E$ equivalence classes of $l$-component string links, resp. links.  
We will also denote by $E$ the map which sends a link or string link to its equivalence class.

Consider the following  set of Axioms for an equivalence relation $E$: 

For $i=1,2$, let $\s_i$ be $l$-component string links with
$E(\s_1)=E(\s_2)$, and let $\Sm_i$ be $2l$-component string
links with $E(\Sm_1)=E(\Sm_2)$. 
\be
\item[(1)]  $E(\hat\s_1)=E(\hat\s_2)$
\item[(2)]  $E(1_l\otimes\s_1)=E(1_l\otimes\s_2)$
\item[(3)]  $E(\s_1 \Sm_1)=E(\s_2\Sm_2)$  
\item[(4)]  $E(\Sm_1 \s_1)=E(\Sm_2\s_2).$  
\item[(5)]  For all string links $\s$, there
is a string link $\s_1$, such that $E(\s\s_1)=E(1_l)$.
\item[(6)]  If $E(L)=E(L')$, then there is an $m$ and a sequence of string links
$\s_i$, for $i=1,\dots, m$, such that $L$ is isotopic to $\hat\s_1$,
and $L'$ is isotopic to $\hat\s_m$, and for all $i$, $1\le i< m$, either
$E(\s_i)=E(\s_{i+1})$, or $\hat\s_i$ is isotopic
to $\hat\s_{i+1}$ (i.e., the equivalence relation $E$ on links is
generated by the equivalence relation of isotopy on links and the
equivalence relation $E$ on string links).
\item[($5'$)]  For all string links $\s$, $E(\s\overline\s)=E(1_l)$.  Here the string link $\overline \s$ is defined by, $\overline \s=
R_t\circ\s \circ R_s$, where $R_s$ and $R_t$ are the reflection
mappings at the source and target.

\ee

\begin{defi}
An equivalence relation satisfying Axioms $(1)-(4)$ is called local.
\end{defi}

We have the following result.
\begin{prop}\label{2.1}
Let $E$ be a local equivalence relation.  For $i=1,2$, let $\s_i$ and  $\s'_i$ be $l$-component string links with
$E(\s_1)=E(\s_2)$ and $E(\s'_1)=E(\s'_2)$, and let $\Sm_i$ be $2l$-component string
links with $E(\Sm_1)=E(\Sm_2)$.   Then $E(\s_1\s'_1)=E(\s_2\s'_2)$ and $E(\s_1 \Sm_1 \s'_1)=E(\s_2\Sm_2\s'_2)$.

The monoidal structures, the left (resp. right) action and the closure mapping all
pass to maps of equivalence classes.  Let $ES^R(l)$ (resp. $ES^L(l)$), denote the right (resp. left) stabilizer of the unit element of $ESL(l)$.  Then $ES^R(l)$ (resp. $ES^L(l)$)
is a submonoid of $ESL(2l)$.  Furthermore, the closure mapping of $ESL(l)$ to $EL(l)$, passes to the set of orbits of the $ES^R(l)$ (resp. $ES^L(l)$) action, i.e., we have a map

$$ {{ESL(l)}\over{S^R(l)}} \longrightarrow  EL(l).$$

If in addition, Axiom $(5)$ holds, 
then the monoid $ESL(l)$ is a group, and $ES^R(l)$ (resp. $ES^L(l)$)
is a subgroup of $ESL(2l)$.  If Axiom $(5')$ holds, then $ES^R(l)=ES^L(l)$.
\end{prop} 

\begin{proof}
By Axiom $(2)$, $E(1_l\otimes\s'_1)=E(1_l\otimes\s'_2)$.  Using Axiom $(3)$, we have that
$E(\s_1\s'_1)=E(\s_1(1_l\otimes\s'_1))=E(\s_2(1_l\otimes\s'_2))=E(\s_2\s'_2)$.

One defines $E(\s_1) E(\s_2)=
E(\s_1\s_2)$.  This is well defined by the above, and the element
$E(1_l)$ is a unit.  One also defines $E(\s) E(\Sm)= E(\s\Sm)$ and 
$E(\Sm) E(\s)= E(\Sm\s)$.
These are well defined, by Axioms $(3)$ and $(4)$, 
and are monoidal actions (of sets).

Suppose $\s$ and $\s'$ define the same element of
$ {{ESL(l)} \over {S^R(l)}}$, i.e., there is $E(\Sm)\in S^R(l)$ such that 
$E(\Sm)E(\s)=E(\s')$.  One has
$E(\hat {\s}')=
E(1_l\Sm\s)=E(1_l 1_{2l}\s)=E(\hat\s)$.  

If Axiom $(5)$ holds, each element $E(\s)$
has a right inverse.  Hence $ESL(l)$ is a group, for all $l$.  

If $E(\Sm)$ belongs to $ES^R(l)$, then $E(\overline\Sm)$ belongs to $ES^L(l)$.  But if Axiom $(5')$ holds, then
$E(\Sm)$ is the inverse of $E(\overline\Sm)$, so also belongs to $ES^L(l)$.  This proves one inclusion and the other is proven similarly.
\end{proof}

\begin{thm}[Structure Theorem for $E$-equivalence]  \hfill
\be 
\item Let $E$ be a local equivalence relation satisfying Axiom $(5)$. 
Then the quotient map 
 \[SL(l)\longrightarrow {{ESL(l)} \over {ES^R(l)}} \]
factors through the closure mapping, i.e., we have a link invariant
 \[ \widetilde{E}: L(l) \longrightarrow {{ESL(l)} \over {ES^R(l)}}.  \]
such that the composite map to $EL(l)$ is $E$.

\item Furthermore, if Axiom $(6)$ also holds, then we have a bijection
$$ {{ESL(l)} \over {S^R(l)}} = EL(l).$$
\ee 
\end{thm}
 
\begin{proof}
 
Suppose $\hat\s$ is isotopic to $\hat\s'$.  By 
Proposition \ref{2.1}, one has that, for some 
$\Sm_0$, $\s'$ is isotopic to $\Sm_0 1_l$ and 
$\s$ is isotopic to $1_l\Sm_0 $. Set $\Sm=
\Sm_0(1_l\otimes \s_1)$, where $E(\s_1)$
satisfies $E(\s\s_1)=E(1_l)$ (and hence also
$E(\s_1\s)=E(1_l)$).  One has that $E(1_l)E(\Sm)=E(1_l \Sm)=
E(\s\s_1)=E(1_l)$, so  
$E(\Sm)\in S^R(l)$.  See Figure \ref{proof1}.  
\begin{figure}[!h]
\includegraphics{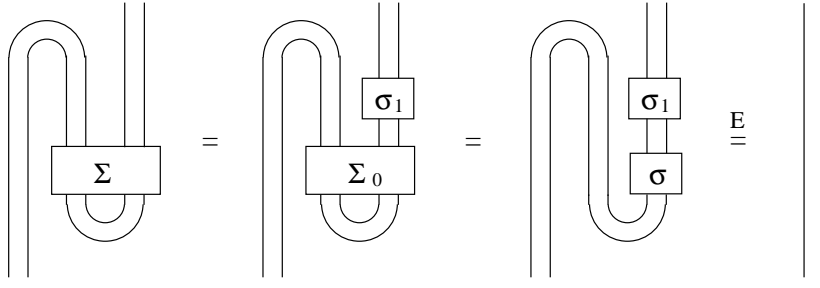} 
\caption{Proof that $\Sigma$ lies in $S^R(l)$. } \label{proof1}
\end{figure}

\noindent Finally, one has that $E(\Sm)E(\s)=E(\Sm\s)=E(\Sm_0)E(\s_1\s)=E(\Sm_0)E(1_l)=
E(\Sm_0 1_l)=E(\s')$.   See Figure \ref{proof2}.  
\begin{figure}[!h]
\includegraphics{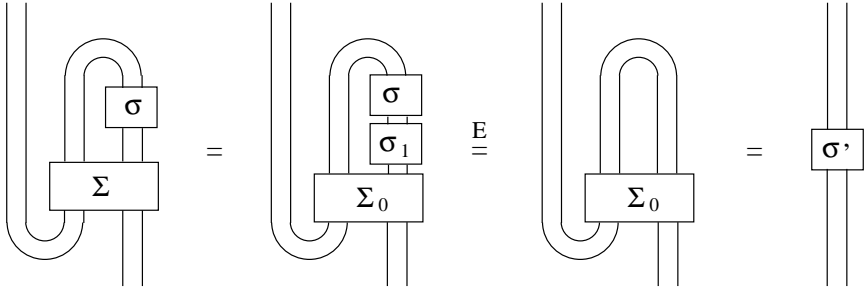} 
\caption{Proof that $\Sigma \sigma$ is equivalent to $\sigma'$.}\label{proof2}
\end{figure}
\noindent This completes the proof of (1).

To see (2), note that we have already shown that if the closures of two string links are isotopic, then they define the same element
of $ {{ESL(l)} \over {S^R(l)}}$.  Thus we have that for all $i$ in Axiom $(6)$, $E(\s_i)$ and $E(\s_{i+1})$ both agree in
$ {{ESL(l)} \over {S^R(l)}}$.  Thus the surjective map from $ {{ESL(l)} \over {S^R(l)}}$ to $ EL(l)$ is injective.

\end{proof}

\section{Structure Theorems for $C_n$-equivalence and for Self-$C_n$-equivalence.} 

We will denote by $FT_n$, the equivalence relation on tangles determined by finite type equivalence up to degree $n$, i.e., $FT_n$-equivalent tangles differ
by an element in the $n+1$st term of the Vassiliev filtration.
In \cite{H}, K. Habiro showed that, \emph{for knots}, $FT_n$-equivalence agrees with another equivalence relation, called $C_{n+1}$-equivalence.  
Habiro conjectured  in \cite{H} that for string links, $FT_n$ equivalence is equivalent to $C_{n+1}$-equivalence. 

Habiro also showed that for links, the result does not hold.  
Note that, since the structure theorem holds for $C_{n+1}$-equivalence, if the equivalence relations were the same both for string links and for links, it would also hold for $FT_n$ equivalence.  However, for $FT_n$ equivalence, the structure theorem does not hold (see Theorem 5.1 and the Borromean ring example of Section 5).  

By definition, two tangles are said to be $C_n$-equivalent, if there is a finite sequence of tree clasper surgeries, of degree greater than or equal to $n$, taking one tangle to the other, up to isotopy.  See \cite{H} for the definition.  (Note that in \cite{H},
a tree clasper is called an admissible, strict tree clasper.)  Here the leaves of the tree can be assumed to be trivial and intersect the tangle in a single point.  
It is known that $C_{n+1}$-equivalent tangles are $FT_n$-equivalent (see \cite[\S 6]{H}).  

By definition, two tangles are said to be Self-$C_n$-equivalent, if there is a finite sequence of tree clasper surgeries, of degree greater than or equal to $n$, taking one tangle to the other, up to isotopy, such that the leaves of each tree are restricted to all intersect the same tangle component.  

\begin{rem} 
Self-$C_n$-Equivalence, for $n=1$, is link-homotopy.  For $n=2$ it is also known as Self-Delta equivalence.  
\end{rem}

$C_n$-equivalence and Self-$C_n$-equivalence are obviously local, i.e., they satisfy Axioms $(1)-(4)$ of Section 3.
Axiom 5 was shown in \cite[Theorem 5.4]{H} for $C_n$-equivalence.

\begin{prop}\label{5.1}
Self-$C_n$-equivalence 
satisfies Axiom $(5)$ of Section 3.
\end{prop}

\begin{prop}\label{5.2}
$C_n$-equivalence and Self-$C_n$-equivalence 
satisfy Axiom $(6)$ of Section 3.
\end{prop}

Applying Theorem 3.2, one has the following result.

\begin{thm}[Structure Theorem for $C_n$-equivalence and Self-$C_n$-equivalence]  
$$ {{C_n SL(l)} \over {C_nS^R(l)}} = C_n L(l).$$
$$ {{\textrm{\textit{Self-C}}_n SL(l)} \over {\textrm{\textit{Self-C}}_nS^R(l)}} = \textrm{\textit{Self-C}}_n L(l).$$
\end{thm}

\begin{proof}[Proof of Proposition \ref{5.2}]
Suppose that $L'$ is obtained from $L$ by surgery on a disjoint union $F$ of tree claspers of degree $\ge n$.  Let $L$ be the closure of $\s$.  Since the disk base for $L$ retracts onto a 1-complex, we may assume it is disjoint from $F$.  Thus $L'$ is the closure of a string link $\s'$, obtained from $\s$ by surgery on a union of tree claspers of degree $\ge n$.  This shows that $C_n$-equivalence (resp. Self-$C_n$-equivalence) for links is implied by $C_n$-equivalence (resp. Self-$C_n$-equivalence) for string links and isotopy.
\end{proof}

Proposition \ref{5.1} is a special case of the following result.

\begin{prop}\label{5.3}
Let $S$ be any surface.
Self-$C_n$-equivalence (and consequently $C_n$-equivalence) classes of string links in $S\times I$ form a group.
\end{prop}

\begin{proof}[Proof of Proposition \ref{5.3}]
The proof is by induction on the number $l$ of components.  For $l=1$, Self-$C_n$-equivalence is $C_n$-equivalence, so we may invoke 
\cite[Theorem 5.4]{H}.  

Suppose the result is true for $l-1$.  Removing the first component from $\s$, 
we have an $l-1$-component string link $\s_0$.  By the induction hypothesis, 
$\s_0$ has an inverse $\s_1$, up to Self-$C_n$-equivalence.  Let $\s'=\s(1_1\otimes\s_1)$.  It suffices to find a right inverse for $\s'$, up to Self-$C_n$-equivalence.
Note that the string link $\s'_0$, obtained from $\s'$ by removing the first component, is Self-$C_n$-equivalent to the trivial string link $1_{l-1}$.  
Thus $1_{l-1}$ is obtainable from 
$\s'_0$ by surgery on a disjoint union $F$ of trees of degree $n$ such that the leaves of each tree are restricted to intersect a single component.

We may assume that $F$ is disjoint from the first component of $\s'$.  Perform surgery on $F$ to obtain from $\s'$ a string link $\s''$.  
As $\s'$ is Self-$C_n$-equivalent to $\s''$, it remains to find a right inverse for $\s''$.  
Note that after removing the first component of $\s''$, we obtain the trivial $l-1$-component string link.  Thus if we remove from $\s''$ the last $l-1$ components, we have a one-component string link $\s''_0$ in $S'\times I$, where $S'$ is the surface obtained from $S$ by removing $l-1$ points.  Since, by the result for $l=1$, the string link $\s''_0$ has a right inverse, up to Self-$C_n$-equivalence, so does $\s''$.

\end{proof}

\section{The Indeterminacy of Finite Type Invariants.} 

In this section we assume the reader is familiar with the notion of finite type invariant as well as the Kontsevich Integral, which is the Universal Finite Type Invariant.  Recall from the last section that we have denoted by $FT_n$, the equivalence relation on tangles determined by finite type equivalence up to degree $n$, i.e., $FT_n$-equivalent tangles differ
by an element in the $n+1$st term of the Vassiliev filtration.

Let us begin with a disturbing fact about finite type invariants of links.  
The Borromean Rings are distinguished from the unlink by the triple Milnor Invariant.  
Unfortunately, this invariant, which is really only defined as an integer when the 
linking numbers of the 2-component sublinks vanish, dies in the space of trivalent Feynman diagrams (also known as Jacobi diagrams)
on 3 circles.  This is because, when passing from 3 intervals to 3 circles, invariants 
of linear combinations of string links, which die upon closure, must also die upon closure for other linear combinations which are equivalent.  

This can be seen using the Kontsevich Integral.  Specifically, in the space of Jacobi diagrams on $3$ intervals we have 
 \[ \includegraphics{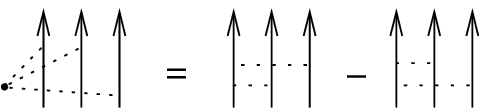} \] 
where the right-hand side is obviously mapped to zero when closing.  
(Recall that the coefficient of the $Y$-shaped diagram on the left-hand side corresponds to the triple Milnor Invariant.)  

To see how this comes about more geometrically, consider the free group on 2 generators as a subgroup of the 3 component pure braid group.
The word $xyx^{-1}y^{-1}$  represents the Borromean rings, after closure.
Since $xyx^{-1}$ and $y$ are conjugate, and thus agree after closure, we see that the quantity 
$xyx^{-1}y^{-1}-1$, which (say after applying the Magnus expansion) is in degree 2 before closure, lies in degree 3 after closure, since it agrees after closure with
the quantity $(xyx^{-1}y^{-1}-1)(y-1)$, which is in degree 3.  (The degree considerations here are valid in the Vassiliev filtration as well).
Thus we see that we can no longer distinguish the Borromean rings from the unlink!

In summary, the indeterminacies of higher order invariants due
to the non-vanishing of lower order ones, propagate to destroy what should be 
invariants of links whose lower order invariants vanish.
We are thus led to a problem of refining the indeterminacies in a less algebraic way.  
We are guided by the structure theorem of the last section.
 
Rationally, it is known, see \cite{HM}, that the set rational finite type equivalence classes of $l$-component string links is a finitely generated torsion free nilpotent group.  Over the integers, it follows from the last section, since $C_{n+1}SL(l)$ is a group and surjects to  $FT_nSL(l)$, that  $FT_nSL(l)$ is also a group.

The set $FT_nSL(2l)$ acts on $FT_nSL(l)$ on the left and right.  Let $FT_nS^R(l)$ denote the stabilizer of $FT_n(1_l)$ under the right action.  
$FT_nL(l)$ denotes the set of $FT_n$ equivalence classes of $l$-component links.

The main result of this paper is the following.
%
\begin{thm}[Structure Theorem for Finite Type Equivalence]  \hfill
\be 
\item The projection mapping,
of $SL(l)$ to the set $\frac{FT_nSL(l)}{FT_nS^R(l)}$ of left $FT_nS^R(l)$ orbits, 
factors through $L(l)$ and thus gives a well defined invariant of links
 \[ \widetilde{FT_n}: L(l) \longrightarrow {{FT_nSL(l)} \over {FT_nS^R(l)}}.  \]
\item The above link invariant lifts the indeterminacies given by finite type invariants of links, i.e., if two links determine the same element of $\frac{FT_nSL(l)}{FT_nS^R(l)}$, then they have the same finite type invariants up to degree $n$.  That is, the above map, $\widetilde{FT_n}$, factors through a (surjective, but not generally injective) map, 
 \[ {{FT_nSL(l)} \over {FT_nS^R(l)}} \longrightarrow  FT_n L(l),  \]
and the composite mapping is 
 \[{FT_n}: L(l) \longrightarrow FT_n L(l).  \]
\ee
\end{thm}
  
\begin{proof}
Axioms $(1)-(4)$ follow from the local definition of the Vassiliev filtration.  Axiom $(5)$ follows from the remark above that $FT_nSL(l)$ is a group.

\end{proof}

\begin{rem}
The analogous theorem also holds if one restricts to the equivalence relation $FT_n^Q$, given by rational invariants of finite type of degree up to $n$.  One can use the local property of the Kontsevich Integral (and the result cited above from \cite{HM}) to give an alternative proof of the Axioms $(1)-(5)$ in this case.

Let $A_{\le n}(l)$ denote the algebra of Jacobi diagrams on $l$ strands of degree up to $n$. The action of
$FT_n^QSL(2l)$ on the set $FT_n^QSL(l)$ is induced, via the Kontsevich Integral, 
by an analogously defined action of $A_{\le n}(2l)$ on $A_{\le n}(l)$, given purely diagrammatically.  (In the definition
of the action of string links, just replace the string links with diagrams.)
Let $A_{\le n}(2l)_1$ be the stabilizer of the unit element in  $A_{\le n}(l)$.  
The stabilizer $A_{\le n}(2l)_1$ contains $FT_n^QS^R(l)$.  It is easily seen that there are surjective maps of the
space of covariants $A_{\le n}(l)/FT_n^QS^R(l)$ to the space of covariants 
$A_{\le n}(l)/A_{\le n}(2l)_1$, and from $A_{\le n}(l)/A_{\le n}(2l)_1$
to the space $A_{\le n}(\coprod_{i=1}^{l} {\bf S}^1_i )$
of diagrams on $l$ circles, up to degree $n$.  Using the link invariance of our theorem, and the 
universal property of the Kontsevich Integral, one can check that these maps are both isomorphisms.  
(We do not have a diagrammatical proof of this fact.)  It follows that one should not 
pass to covariants to try to refine finite type invariants of links!
\end{rem}

We conclude this section with several problems.   

\begin{Problem}
Use the `unipotent' action to write an algorithm, analogous to Xiao-Song's link-homotopy algorithm, `calculating' whether or not two (string) links determine the same element in the orbit space.
\end{Problem}
\begin{Problem} 
Does the full Kontsevich Integral for links (or integrally, modulo the intersection of the Vassiliev filtration) `recapture' the information lost at each finite level? 
(For example, the triple Milnor Invariant dies, but its cube does not. But of course the degree is now 6 and not 2.)
\end{Problem}

\end{document}